\documentclass[12pt,reqno,psfig,openbib]{amsart}

\usepackage[pagebackref=true, colorlinks=true, citecolor=blue]{hyperref}
\usepackage{amssymb,amsmath,graphicx,amsfonts,euscript}
\usepackage{fancyhdr}
\usepackage{epsfig}
\usepackage{setspace}
\usepackage{dsfont}
\usepackage[all]{xy}
\usepackage{amsmath}
\usepackage{enumitem}

\usepackage{caption}
\usepackage{subcaption}

\usepackage{soul}

\usepackage{float}

 \usepackage{pst-grad} 
 \usepackage{pst-plot} 

\numberwithin{equation}{section}
\numberwithin{figure}{section}

 \usepackage[left=1.05in,top=1.15in,right=1.05in]{geometry}
\usepackage{verbatim}

\newcommand{\spacer}{\vspace{1.5mm}}

\date{\today}

\renewcommand{\epsilon}{\varepsilon}
\newcommand{\ep}{\varepsilon}

\newcommand{\pri}{\prime}
\newcommand{\sig}{\sigma}
\newcommand{\kap}{\kappa}
\newcommand{\blds}{\boldsymbol}

\begin{document}

\title[Semi-analytic PINN methods for singular perturbations]
{Semi-analytic PINN methods for singularly perturbed boundary value problems}

\author[G.-M. Gie, Y. Hong, and C.-Y.  Jung]
{Gung-Min Gie$^{1}$, 
Youngjoon Hong$^{2}$, and 
Chang-Yeol  Jung$^3$}
\address{$^1$ Department of Mathematics, University of Louisville, Louisville, KY 40292}
\address{$^2$ Department of Mathematics, Sungkyunkwan University}
\address{$^3$ Department of Mathematical Sciences, Ulsan National Institute of Science and Technology, 
Ulsan 44919, Korea}
\email{gungmin.gie@louisville.edu}
\email{hongyj@skku.edu}
\email{cjung@unist.ac.kr}

\begin{abstract}
We propose a new semi-analytic physics informed neural network (PINN) to
solve singularly perturbed boundary value problems. The PINN is a scientific machine
learning framework that offers a promising perspective for finding numerical solutions to
partial differential equations. The PINNs have shown impressive performance in solving
various differential equations including time-dependent and multi-dimensional equations
involved in a complex geometry of the domain. However, when considering stiff differential equations, neural networks in general fail to capture the sharp transition of solutions,
due to the spectral bias. To resolve this issue, here we develop the semi-analytic PINN
methods, enriched by using the so-called corrector functions obtained from the boundary
layer analysis. Our new enriched PINNs accurately predict numerical solutions to the singular perturbation problems. Numerical experiments include various types of singularly
perturbed linear and nonlinear differential equations.
\end{abstract}

\maketitle

\tableofcontents

\section{Introduction}\label{S. Intro}
 
Neural networks have been widely studied and used 
for approximating 
solutions to differential equations; 
see, e.g.,
\cite{
BN18,
CX21, 
CCWX22, 
EY18, 
KZK21, 
KDJH21, 
RK18, 
RPK19, 
RBPK17, 
QCJX21}.
In this research direction, 
many (unsupervised) neural networks without training datasets have been successfully 
developed, e.g.,  
physics informed neural networks (PINNs)  
\cite{
BE21, 
KZK19, 
KDJH21, 
HL21}, 
deep Ritz method (DRM)  
\cite{Yu18}, 
and deep Galerkin method (DGM)  
\cite{SS18} 
where 
the loss function is defined by 
using 
a certain residual from the differential equation under consideration.
Especially the 
PINNs use  collocation points in the space-time domain 
as the training data set, 
and hence 
the PINNs are suitable for solving 
time dependent, 
multi-dimensional equations involved in a complex geometry 
of the domain, 
\cite{
KZK19, 
LMMK21, 
LDKRKS20, 
SDBSK20, 
YMK21, 
MZK20, 
KKLPWY21, 
WSP22, 
XZRW21}.

Compared to many other types of neural networks, 
some advantages of using the PINNs in the study of differential equations 
include, but not limited to, 
first the fact that 
they are unsupervised learning processes and 
hence the exact solution of a model differential equation is not a-priori required for the learning process. 
The exact solution is usually used only when we measure 
the error between the exact solution and the approximate solution obtained by the PINNs, and  
hence the PINNs work as like the traditional numerical methods for differential equations, e.g., finite difference, finite elements, and so on. 
Another big advantage of using PINNs 
(over the traditional numerical methods) is their flexibility applied to many different types of differential equations 
mainly because the model equation under consideration 
is used only when the loss is computed.  
However, being unsupervised learning process, a repeat of learning is required for different data, i.e., external force, or initial data, and each learning takes a significant amount of time, compared to the relatively short time cost of most traditional numerical methods. 
\bigskip

In this article, 
we  
approximate the solutions of 
various 1D boundary value problems, especially when the highest order derivative appearing in each equation is multiplied by a 
small parameter $\ep > 0$. 
More precisely, 
we consider the 1D elliptic differential equations in the form,
\begin{equation}\label{e:model_gen}
    F(x;\, 
    u,\, 
    u_x, \,
    \ep u_{xx}
    )
    =
    f(x), \quad
    0 < x < 1,
\end{equation}
supplemented with a Dirichlet boundary condition, 
\begin{equation}\label{e:model_gen_BC}
    u = 0, \quad  
    x=0,1.
\end{equation}

A singularly perturbed boundary value problem, such as our problem (\ref{e:model_gen}) - (\ref{e:model_gen_BC}), 
is well-known to 
generate a thin layer 
near the boundary {\em (called the boundary layer)},   
in which 
a sharp transition of the solution occurs.  
A large literature has been developed on the mathematical theory of singular
perturbations and boundary layers; see, e.g., 
\cite{
Book16, 
Ho95, 
OM08, 
SK87}. 
Concerning 
numerical approximation of the singular perturbation problem,     
a very large computational error is created near the boundary,    
due to the stiffness of solution inside the boundary layer.  
Thus, 
to achieve a sufficiently accurate approximation of the solution near the boundary,  
a massive mesh refinement is usually required, near the boundary, for the 
most classical numerical schemes. 
Instead of introducing massive mesh
refinements, 
new semi-analytic methods have been proposed, see, e.g., 
\cite{CJL19, 
GJL1, 
GJL2, 
HK82, 
HJL13, 
HJT14}.  
The
main component of this semi-analytic method is 
{\em enriching } the basis of traditional numerical methods, e.g. finite elements, finite volumes and so on, 
by 
adding a global basis function, 
called the 
{\em corrector},  
which 
describes the singular behavior of the solution inside boundary layers. 
Such semi-analytic methods have proven to be highly efficient without any help of mesh refinement near the boundary.
\bigskip

Our main goal  is 
to construct 
a {\em semi-analytic} physics-informed neural networks (PINNs), 
{\em enriched} by using   
the so-called {\em corrector} functions; see Section \ref{S.Sing_eqns} and \ref{S.Sing_eqns_non} below. 
Toward this end, we first 
briefly recall the PINNs (well-developed in earlier works; see, e.g., 
\cite{LMMK21}) 
for our model equation in (\ref{e:model_gen}) - (\ref{e:model_gen_BC}):

For  
an $L$-layer Neural Network (NN) 
(or $L-1$-hidden layer NN),   
the $\ell$-th layer with $N_{\ell}$ neurons is  denoted by  $\mathcal{N}^{\ell}, \, 0 \leq \ell \leq L$. 
Writing the weight matrix and the bias vector at each $\ell$-th layer as $\blds{W}^{\ell} \in \mathbb{R}^{N_{\ell} \times N_{\ell-1}}$, $1 \leq \ell \leq L$,  
$\blds{b}^{\ell} \in \mathbb{R}^{N_{\ell}}$,  $1 \leq \ell \leq L$, 
and 
${b}^{L} \in \mathbb{R}$
(with $N_0 = N_{L} = 1$ for our 1D problem (\ref{e:model_gen}) - (\ref{e:model_gen_BC})),  
we use the feed- forward neural network (FNN) with an activation function $\sigma$, and 
we recursively define  
\begin{equation}\label{e:PINN_gen}
\begin{array}{rl}
      {input\,\,\, layer:}
            &
            \mathcal{N}^0({ x}) = { x} \in \mathbb{R},\\
      {hidden\,\,\, layers:}
            &
            \mathcal{N}^\ell ({  x})
            = \sigma( {\blds W}^{\ell}         \mathcal{N}^{\ell - 1} ({x})  + {\blds b}^{\ell}) 
            \in \mathbb{R}^{N_{\ell}}, 
            \quad 1 \leq \ell \leq L-1,\\
      {output\,\,\, layer:}
            &
            \mathcal{N}^L({ x})
            =
            {\blds W}^{L}\mathcal{N}^{L - 1} ({ x})  + {b}^{L}\in \mathbb{R}.
\end{array}
\end{equation}

Using the NN above, 
we construct an approximation $\hat{u} = \hat{u}( {x}; \, \blds{\theta})$
where the parameters 
$\blds{\theta} = \{ \blds{W}^{\ell}, \, \blds{b}^{\ell} \}_{1 \leq \ell \leq L}$ is the set of the weight matrices and bias vectors in the neural network. 
In the next step, the NN approximation $\hat{u}$ 
is restricted to satisfy the constraint  imposed 
by the PDE and boundary conditions (\ref{e:model_gen}) - (\ref{e:model_gen_BC}). For this purpose, 
we prepare the training set $\mathcal{T} \subset [0, 1]$ which consists of two sets of scattered points, $\mathcal{T}_I 
\subset (0, 1)$ and $\mathcal{T}_B \subset \{ 0, \, 1\}$, 
and define the corresponding loss function, using the weighted $L^2$ norm, by 
\begin{equation}\label{e:loss_PINNs_gen}
    \mathcal{L}({\blds \theta}, \, \mathcal{T})
        =
            w_I \mathcal{L}_I({\blds \theta}; \, \mathcal{T}_I)
            +
            w_B \mathcal{L}_B({\blds \theta}; \, \mathcal{T}_B), 
\end{equation}
where 
\begin{equation}\label{e:loss_PINNs_gen_FB}
\begin{array}{l}
\displaystyle
    \mathcal{L}_I({\blds \theta}, \, \mathcal{T}_I)
        =
            \frac{1}{|\mathcal{T}_I|}
            \sum_{ { x}  \in \mathcal{T}_I}
            \|
            F(x;\,  \hat{u},\,  \hat{u}_x, \, \ep \hat{u}_{xx}
            ) -    f 
            \|_{L^2{(0, 1)}}^2,\\
\displaystyle
    \mathcal{L}_B({ \blds \theta}, \, \mathcal{T}_B)
        =
            \frac{1}{|\mathcal{T}_B|}
            \sum_{{ x} \in \mathcal{T}_B}
            \|
            \hat{u}({ x})
            \|_{l^2{\{ 0, 1 \}}}^2,
\end{array}
\end{equation}
and $w_F$ and $w_F$ are the certain weight parameters. 
The loss in (\ref{e:loss_PINNs_gen}) indeed involves the derivatives of $\hat{u}$ and it is handled via the so-called automatic differentiation (AD). 
In the last step, the procedure of searching for a good $\blds \theta$  by minimizing the loss $\mathcal{L}({\blds \theta}, \, \mathcal{T})$ is called ``training" 
where we usually use a gradient-based optimizer 
such as gradient descent, Adam, or L-BFGS.
\bigskip

Our approach to construct a {\it two-layer} NN solving (\ref{e:model_gen}) - (\ref{e:model_gen_BC}) 
is closely related to the PINNs above, but it is a bit different as explained below;   
see an earlier work \cite{LLF98} where 
the NN closely related to ours is introduced:  

To obtain an approximate solution $\overline{u}$ to (\ref{e:model_gen}) - (\ref{e:model_gen_BC}), 
we employ a simple 2-layer NN (denoted by $\hat{u}$) multiplied by $x(x-1)$ to enforce the boundary condition (\ref{e:model_gen_BC}), 
in the form:
\begin{equation}\label{e:app_sol_mod}
    \overline{u}(x ; \, {\blds \theta}) 
        =
            x(x-1) \, 
            \hat{u}({ {x}}; \, {\blds \theta}), 
\end{equation}
where $\hat{u}$ is defined by the 2-layer NN, 
\begin{equation}\label{e:PINN_modified}
\begin{array}{rl}
      {input\,\,\, layer:}
            &
            \mathcal{N}^0({ x}) = { x} \in \mathbb{R},\\
      {hidden\,\,\, layer :}
            &
            \displaystyle 
            \mathcal{N}^1 ({ x})
            = 
                \sigma( {\blds W}^1 {  x}  + {\blds b}^1 ) 
            \in \mathbb{R}^{N_1},\\
      {output\,\,\, layer:}
            &
            \mathcal{N}^2({  x})
            =
            {\blds W}^2  \mathcal{N}^{1} ({ x})  \in \mathbb{R}, 
\end{array}
\end{equation}
i.e., 
\begin{equation}\label{e:NN_modified}
    \hat{u}({{x}}; \, {\blds \theta})
    =
          {\blds W}^2
          \sigma( {\blds W}^1 {x}  + {\blds b}^1), 
\end{equation}
for the parameters 
${\blds \theta} = \{  
                    {\blds W}^1 
                        = (W^1_1 \, \cdots \, W^1_{N_1})^T, \,
                    {\blds W}^2
                        = (W^2_1 \, \cdots \, W^2_{N_1}), \, 
                    {\blds b}^1
                        = (b^1_1 \, \cdots \, b^1_{N_1})^T\}
    \in \mathbb{R}^{N_1 \times 1}
        \times
        \mathbb{R}^{1 \times N_1}
        \times
        \mathbb{R}^{N_1}.
    $
Here we choose the activation function as 
the logistic sigmoid, 
\begin{equation}\label{e:sigmoid}
    \sigma(x) 
        =
            1/(1 + e^{-x}). 
\end{equation}

Because the boundary condition   (\ref{e:model_gen_BC}) is already embedded in 
the approximate solution (\ref{e:app_sol_mod}), we define the loss function as 
\begin{equation}\label{e:loss_PINNs_mod}
    \mathcal{L}({\blds \theta}, \, \mathcal{T})
        =
            \frac{1}{|\mathcal{T} |}
            \sum_{ {\blds x}  \in \mathcal{T} }
            \|
            F(x;\,  \overline{u},\,  \overline{u}_x, \, \ep \overline{u}_{xx}
            ) -    f 
            \|_{L^2(0, 1)}^2, 
\end{equation}
where 
the training set $\mathcal{T}$ is chosen as a set of scattered points in 
{$(0, 1)$}. 

One big difference between the usual PINNs and 
our 2 layer modified PINNs is on 
the computing the derivatives of the NN approximation in the loss function. 
Using a fact on the sigmoid function that 
\begin{equation}\label{e:derivative_sigmoid}
\begin{array}{rl}
    \sig^\pri(x)
        &
        \spacer
        =
            \sig(x) \big(
                        1 - \sig(x)
                    \big),\\
    \sig^{\pri \pri}(x)
        &
        =
            \sig (x) 
                    \big(
                        1 - \sig(x)
                    \big)
                     \big(
                        1 - 2 \sig(x)
                    \big), 
\end{array}
\end{equation}
we compute
\begin{equation}\label{e:NN_modified_der}
\begin{array}{rl}
    \dfrac{d}{dx}\hat{u}({{x}}; \, {\blds \theta})
    &
    \displaystyle 
    \spacer
    =
        \sum_{j=1}^{N_1}
          W^2_j
          W^1_j
          \sigma^\pri( W^1_j {x}  + b^1_j)\\
    &
    \displaystyle \spacer
    =
        \sum_{j=1}^{N_1}
          W^2_j
          W^1_j
          \sig ( W^1_j {x}  + b^1_j)
          \big(
                1 - \sig(W^1_j {x}  + b^1_j)
          \big),\\
\dfrac{d^2}{dx^2}\hat{u}({{x}}; \, {\blds \theta})
    &
    \displaystyle 
    \spacer
    =
        \sum_{j=1}^{N_1}
          W^2_j
          (W^1_j)^2
          \sig ( W^1_j {x}  + b^1_j)
          \big(
                1 - \sig(W^1_j {x}  + b^1_j)
          \big)
          \big(
                1 - 2\sig(W^1_j {x}  + b^1_j)
          \big).
\end{array}
\end{equation}
Thanks to the simple structure (\ref{e:app_sol_mod}) of our 2 layer modified PINNs,  
we use 
the derivatives of $\hat{u}({{x}}; \, {\blds \theta})$ above, 
and compute explicitly  
the loss function (\ref{e:loss_PINNs_mod})  
without using the automatic differentiation (AD); 
hence our new 2 layer modified PINNs do not rely on the AD and 
their computational errors are independent of the AD.  
In addition, see the sections below where we use (\ref{e:app_sol_mod}) 
(or the proper modifications of (\ref{e:app_sol_mod}) by using the so-called boundary layer correctors)  
to define the loss function of each example we consider in this article.
\bigskip

In general, when the target function contains high-frequency components, 
i.e., when $\ep>0$ is small in our model probelem (\ref{e:model_gen}), the PINN algorithms often fail to converge to the desirable solutions, 
because the so-called spectral bias phenomenon \cite{CF21, RB19}. 
Since the general learning process of neural networks rely on the smooth prior, the spectral bias leads to a failure to capture sharp transitions accurately or singular behaviors of the target solution function.
More precisely, while the neural networks tend to learn low frequency components, it requires much time to fit high frequency components.  
Hence, without care, neural networks cannot fit the sharp transition caused by the boundary layer.
In our algorithm, we split our test function into two parts, slow and fast gradient components, and learn effectively the two components simultaneously.
This approach is motivated from an enriched basis method in numerical PDEs \cite{GJL1, GJL2, HJL13, HJT14}.
\bigskip

In Section \ref{S.reg} below, 
we first verify that 
the 2 layer modified PINNs work well for 
the (regular) problem (\ref{e:model_gen}) - (\ref{e:model_gen_BC}) when 
the parameter is not small, e.g., $\ep =1$. 
Then, in Section \ref{S.Sing_eqns} and \ref{S.Sing_eqns_non}, 
as the main work of our article, 
we consider different types of linear and non-linear singular perturbation problems in the form of 
(\ref{e:model_gen}) - (\ref{e:model_gen_BC}) 
when the parameter $\ep$ is small. 
As we will see below, 
the 2 layer modified PINNs (as well as the usual PINNs) 
do not capture well the singular behavior of perturbation problems, 
caused by 
the so-called boundary layers, and hence they 
fail to produce an accurate approximate solution for each 
singular perturbation.  
To overcome this difficulty of boundary layers, 
we first perform the boundary layer analysis of each singular perturbation problem, and find 
the so-called corrector function 
that exhibits the singular behavior 
of the problem inside the boundary layer. 
Then we  
construct our new semi-analytic Neural Networks,   
enriched 
by embedding the corrector function in the structure of 2 layer modified PINNs. 
Numerical simulations for each example below confirm that 
our new enriched 2 layer modified PINNs captures naturally the singular behavior of boundary layers and produce 
accurate approximations of 
singularly perturbed boundary value problems considered in this article.

\subsection{Elliptic and hyperbolic equations}\label{S.reg}
Our first task is to numerically confirm that 
the modified PINN (\ref{e:app_sol_mod}) 
works as good as the usual PINN (\ref{e:PINN_gen}) 
for a certain class of boundary value problems. 
To this end, 
{we introduce} the following boundary value problem, 
\begin{equation}\label{e:reg_rea_con_dif}
    \left\{
    \begin{array}{rl}
        - a u_{xx} - b u_x + c  u = f(x), & 0 < x <1, \\
        u = 0, & x=0, 1.
    \end{array}
    \right.
\end{equation}
By setting the coefficients $(a, b, c)$ as  
$(0, -1, 0)$, 
$(1, 1, 0)$, and 
$(1, 0, 1)$, 
we consider (\ref{e:reg_rea_con_dif}) as 
the 
{\it hyperbolic}, {\it convection-diffusion}, and {\it reaction-diffusion} equation.

We apply our 
2-layer modified PINN with the sigmoid function, 
defined in (\ref{e:app_sol_mod}), to approximate solutions to 
the boundary value problem (\ref{e:reg_rea_con_dif}). 
The computation results below 
in Fig. \ref{fig:plain} 
show 
that our 2 layer modified PINNs 
works well and produce accurate approximations for the solutions to (\ref{e:reg_rea_con_dif}); 
see \cite{LLF98, KKLPWY21, LMMK21, RK18} as well for the comparable computational results of the usual PINNs.

\begin{figure}[h!]
     \centering
     \begin{subfigure}[b]{0.3\textwidth}
         \centering
         \includegraphics[width=\textwidth]{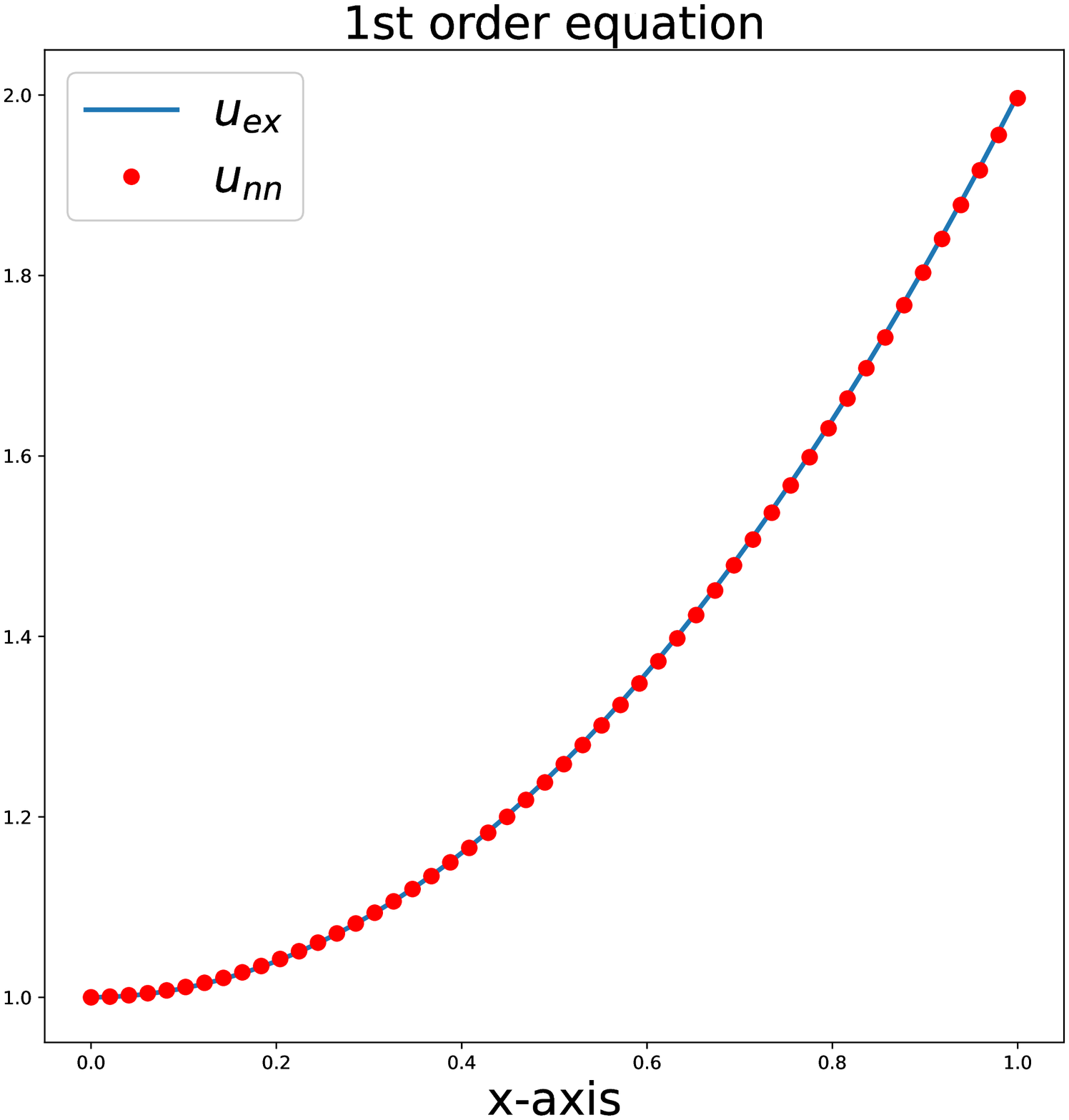}
         \caption{1st order equation with $a=0, c=0, b= -1$}
     \end{subfigure}
     \hfill
     \begin{subfigure}[b]{0.3\textwidth}
         \centering
         \includegraphics[width=\textwidth]{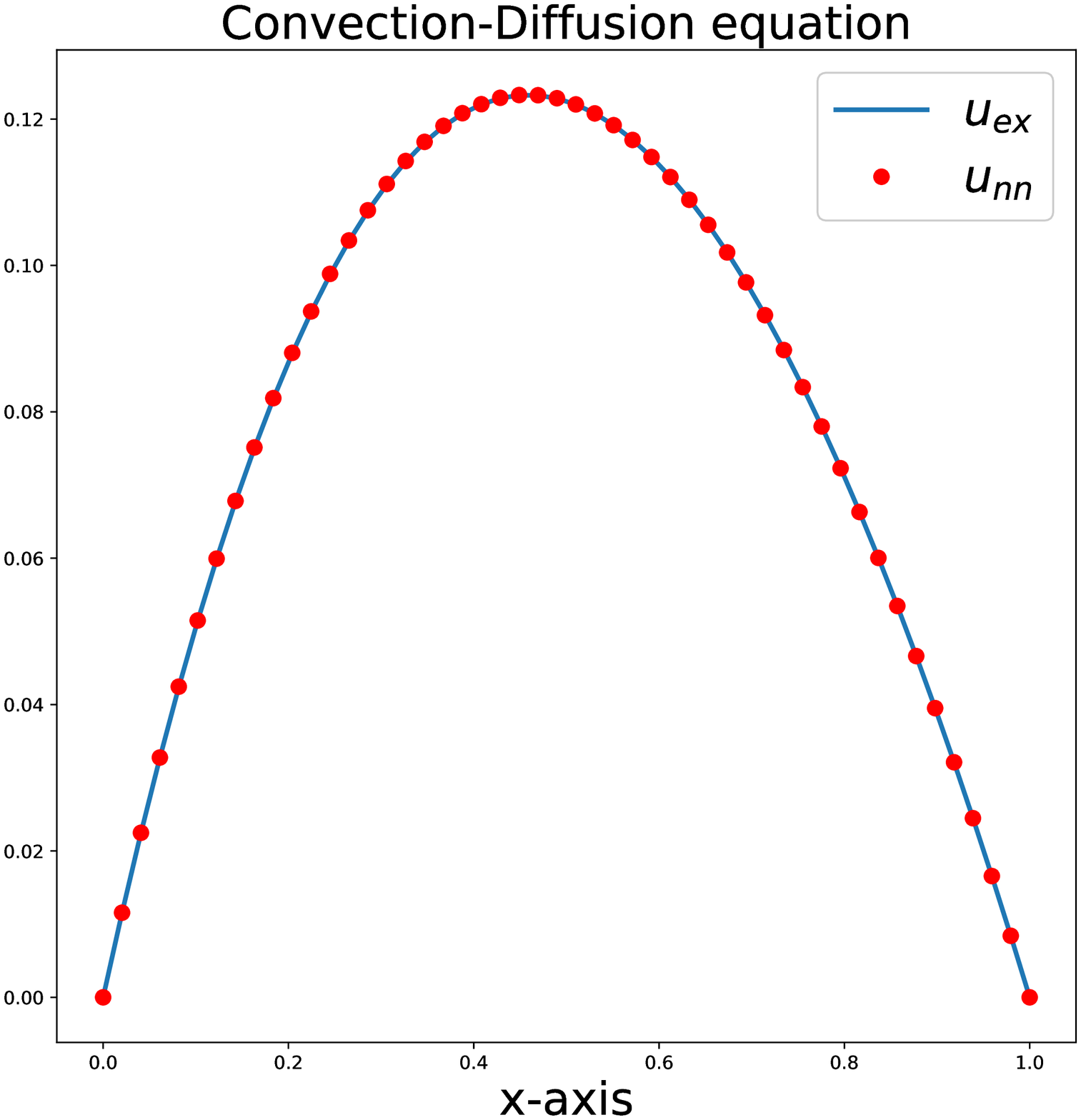}
         \caption{Convection-diffusion eqn with $a=1, b=1,c =0$.}
     \end{subfigure}
     \hfill
     \begin{subfigure}[b]{0.30\textwidth}
         \centering
         \includegraphics[width=\textwidth]{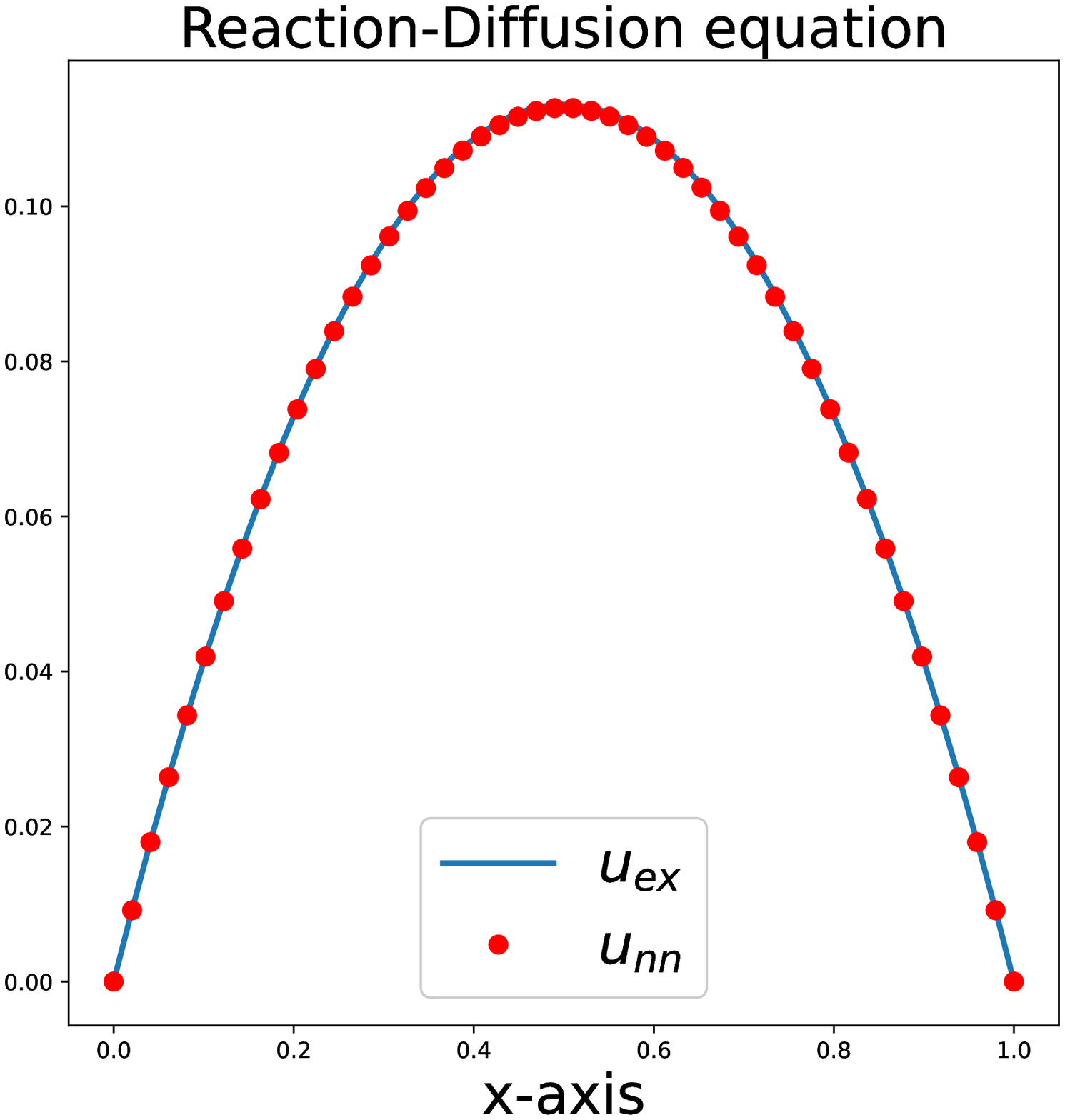}
         \caption{Reaction-diffusion eqn with $a=1, b=0,c =1$.}
     \end{subfigure}
        \caption{Comparison on the exact solution of \eqref{e:reg_rea_con_dif} and our neural network approximation of \eqref{e:reg_rea_con_dif} using the 2-layer modified PINN. 
        For numerical experiments, we set $f=1$ in \eqref{e:reg_rea_con_dif}}
        \label{fig:plain}
\end{figure}

\section{Enriched Neural Network: linear singularly perturbed boundary value problems}\label{S.Sing_eqns}

In this section, we study a linear singularly perturbed boundary value problem. We
first consider the convection-diffusion equations where the boundary layer occurs near the
outflow boundary only. The reaction-diffusion equations are then investigated where the boundary layer takes place at each boundary.

\subsection{Convection-diffusion equations}\label{S.con_dif_eqn}

As a first project, 
we consider the following singularly perturbed convection-diffusion equation, 
\begin{equation}\label{e:con_dif_eqn}
    \left\{
    \begin{array}{rl}
        - \ep u^\ep_{xx} - u^\ep_x  = f(x), & 0 < x <1, \\
        u^\ep = 0, & x=0, 1.
    \end{array}
    \right.
\end{equation}
Our main object here is 
to construct the semi-analytic 
PINN for the singular perturbation problem above, and 
compare its performance with 
the PINN:

In order to find the so-called corrector, 
which is suitable to 
enrich the 2-layer PINNs for the problem (\ref{e:con_dif_eqn}) above, 
we 
first notice, by formally replacing $\ep$ by $0$ in (\ref{e:con_dif_eqn})$_1$, that
the limit problem of (\ref{e:con_dif_eqn}) at the vanishing 
diffusivity at $\ep = 0$ is
\begin{equation}\label{e:con_dif_eqn_limit}
    \left\{
    \begin{array}{rl}
        \spacer
        - u^0_x  = f(x), & 0 < x <1, \\
        u^0 = 0, & x= 1.
    \end{array}
    \right.
\end{equation}
Here we impose the so-called {\it inflow} boundary condition at $x=1$ for $u^0$ (see, e.g., \cite{Book16} for more information), and hence we find the formal limit $u^0$ in the form, 
\begin{equation}\label{e:con_dif_u^0}
    u^0 (x) 
        =
            \int_0^1 f \, dx. 
\end{equation}

Now, by performing the matching asymptotics for (\ref{e:con_dif_eqn}), 
we observe that 
the boundary layer of size $\ep$ occurs near the {\it outflow} boundary at $x = 0$, 
and  
find that  
the asymptotic equation 
(with respect to the small $\ep$) for the corrector 
$\varphi$, 
which approximate the difference $u^\ep - u^0$, 
is given 
in the form,
\begin{equation}\label{e:cor_eqn_con-dif}
    \left\{
        \begin{array}{rl}
            -\ep \varphi_{xx} - \varphi_x = 0, 
            &
            0 < x < 1, \\
            \varphi= - u^0 , 
            &
            x = 0.
        \end{array}
    \right.
\end{equation}

It is well-known, see, e.g., \cite{Book16},  that the corrector $\varphi$ is 
given in the form
\begin{equation}\label{e:cor__con-dif}
    \varphi(x)
        =
            -u^0(0) \, e^{-x/\ep}
            + e.s.t., 
\end{equation}
where the $e.s.t.$ denotes an exponentially small term 
with respect to the small perturbation parameter $\ep$. 

Employing the simple energy estimates on the difference 
$u^\ep - (u^0 + \varphi)$, and 
then using the smallness of the corrector $\varphi$, 
we then 
further notice that 
\begin{equation}\label{e:con__con-dif}
    \begin{array}{l}
    \spacer
        \|
            u^\ep - (u^0 + \varphi)
        \|_{L^2((0, 1))}
            \leq
                \kap \ep, \\
        \|
            u^\ep -  u^0  
        \|_{L^2((0, 1))}
            \leq
                \kap \ep^{\frac{1}{2}}, 
    \end{array}
\end{equation}
for a constant $\kap>0$ independent of $\ep$. 
The convergence results above denote 
first that the diffusive solution $u^\ep$ converges, as $\ep$ tends to $0$, to the limit solution $u^0$ as fast as $\ep$. 
Second, we also infer from the convergence results 
that the corrector $\varphi$ 
exhibits the singular behavior of $u^\ep$ at a small diffusivity $\ep$, i.e., 
the diffusive solution $u^\ep$ is decomposed into the sum of 
fast (decaying part) $\varphi$, 
and the slow part $u^0$.
\bigskip

Inspired by the analysis above, 
we modify the 2-layer PINN  by
using the profile of the corrector $\varphi$, 
and define 
our new {\it semi-analytic 2-layer PINN} for the problem (\ref{e:con_dif_eqn}) in the form,
\begin{equation}\label{e:app_sol_enriched_con-def}
    \widetilde{u}(x ; \, {\blds \theta}) 
        =
            (x-1) \, 
            ( 
                \hat{u}({ {x}}; \, {\blds \theta})
            -
                \hat{u}({ {0}}; \, {\blds \theta}) \,
                e^{-x/\ep}
            ), 
\end{equation}
where $\hat{u}$ is exactly the same as in (\ref{e:NN_modified}), which 
was used for the modified 2-layer PINN  approximation $\overline{u}$ 
in (\ref{e:app_sol_mod}). 
Here we observe that 
the boundary value of $\widetilde{u}$ at $x = 0$ is ensured to be zero,  
thanks to the factor 
$- \hat{u}({ {0}}; \, {\blds \theta})$
multiplied to the exponentially decaying boundary layer function; 
see and compare the difference of 
(\ref{e:app_sol_enriched_con-def}) and 
(\ref{e:app_sol_mod}) 
with or without the factor $x$.

The loss for this problem (\ref{e:con_dif_eqn}) is defined by 
\begin{equation}\label{e:loss_PINNs_mod_con-dif}
    \mathcal{L}({\blds \theta}, \, \mathcal{T})
        =
            \frac{1}{|\mathcal{T} |}
            \sum_{ {\blds x}  \in \mathcal{T} }
            \|
            -\ep \widetilde{u}_{xx}
            - \widetilde{u}_x
            - f 
            \|_{L^2(0, 1)}^2, 
\end{equation}
where 
the training set $\mathcal{T}$ is chosen as a set of scattered points in 
{$(0, 1)$}. 
The derivatives of our enriched 2 layer approximation $\widetilde{u}$ are given by 
\begin{equation}\label{e:2layer_app_enriched_der}
\begin{array}{l}
    \spacer
    \dfrac{d}{dx} \widetilde{u}
        =
            ( 
                \hat{u}({ {x}}; \, {\blds \theta})
            -
                \hat{u}({ {0}}; \, {\blds \theta}) \,
                e^{-x/\ep}
            )
            +
            (x-1)
            \Big( 
                \hat{u}_x({ {x}}; \, {\blds \theta})
            + \dfrac{1}{\ep}
                \hat{u}({ {0}}; \, {\blds \theta}) \,
                e^{-x/\ep}
            \Big),\\
    \dfrac{d^2}{dx^2} \widetilde{u}
        =
            2 \Big( 
                \hat{u}_x({ {x}}; \, {\blds \theta})
            + \dfrac{1}{\ep}
                \hat{u}({ {0}}; \, {\blds \theta}) \,
                e^{-x/\ep}
            \Big)
            +
            (x-1)
            \Big( 
                \hat{u}_{xx}({ {x}}; \, {\blds \theta})
            - \dfrac{1}{\ep^2}
                \hat{u}({ {0}}; \, {\blds \theta}) \,
                e^{-x/\ep}
            \Big).
\end{array}
\end{equation}
Using the fact that 
the exponentially decaying function $e^{-x/\ep}$ 
satisfies the model equation (\ref{e:con_dif_eqn}), 
we find that, for 
the loss function in
(\ref{e:loss_PINNs_mod_con-dif}), 
\begin{equation}\label{e:loss_eqn_con-dif}
\begin{array}{l}
\spacer
    -\ep \widetilde{u}_{xx}
    - \widetilde{u}_x
    - f \\
    \qquad
        =
            - \ep (x-1) \hat{u}_{xx} ({ {x}}; \, {\blds \theta})
            - (x - 1 + 2\ep) \hat{u}_{x} ({ {x}}; \, {\blds \theta})
            - \hat{u} ({ {x}}; \, {\blds \theta})
            -
                \hat{u} (0; \, {\blds \theta}) e^{-x/\ep}
            - f
            ,
\end{array}
\end{equation}
where the derivatives of $\hat{u}$ are given in (\ref{e:NN_modified_der}). 
One important remark here  
is that 
all the terms in (\ref{e:loss_eqn_con-dif}) stay bounded as 
$\ep$ vanishes, 
and hence 
our enriched 2 layer PINN 
produces 
an 
accurate 
approximation $\widetilde{u}$ for (\ref{e:con_dif_eqn}),  
independent of the small parameter $\ep$.
\bigskip

\begin{figure}[H]
    \centering
    \includegraphics[width=0.45
    \textwidth]{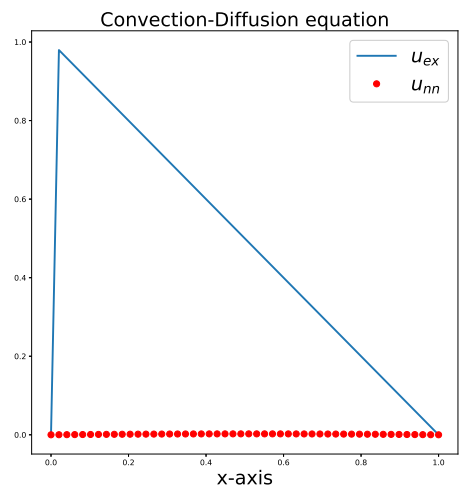}
    \caption{The exact solution and the approximation solution of the convection-diffusion equation \eqref{e:con_dif_eqn} with $\epsilon = 10^{-4}$ are displayed. 
     Here, the (plain) modified 2-layer PINN \eqref{e:app_sol_mod} is used for the predicted solution.
     The predicted solution is not close to the corresponding exact solution.}
    \label{fig:classical_NN}
\end{figure}

\begin{figure}[H]
    \centering
    \includegraphics[width=0.5\textwidth]{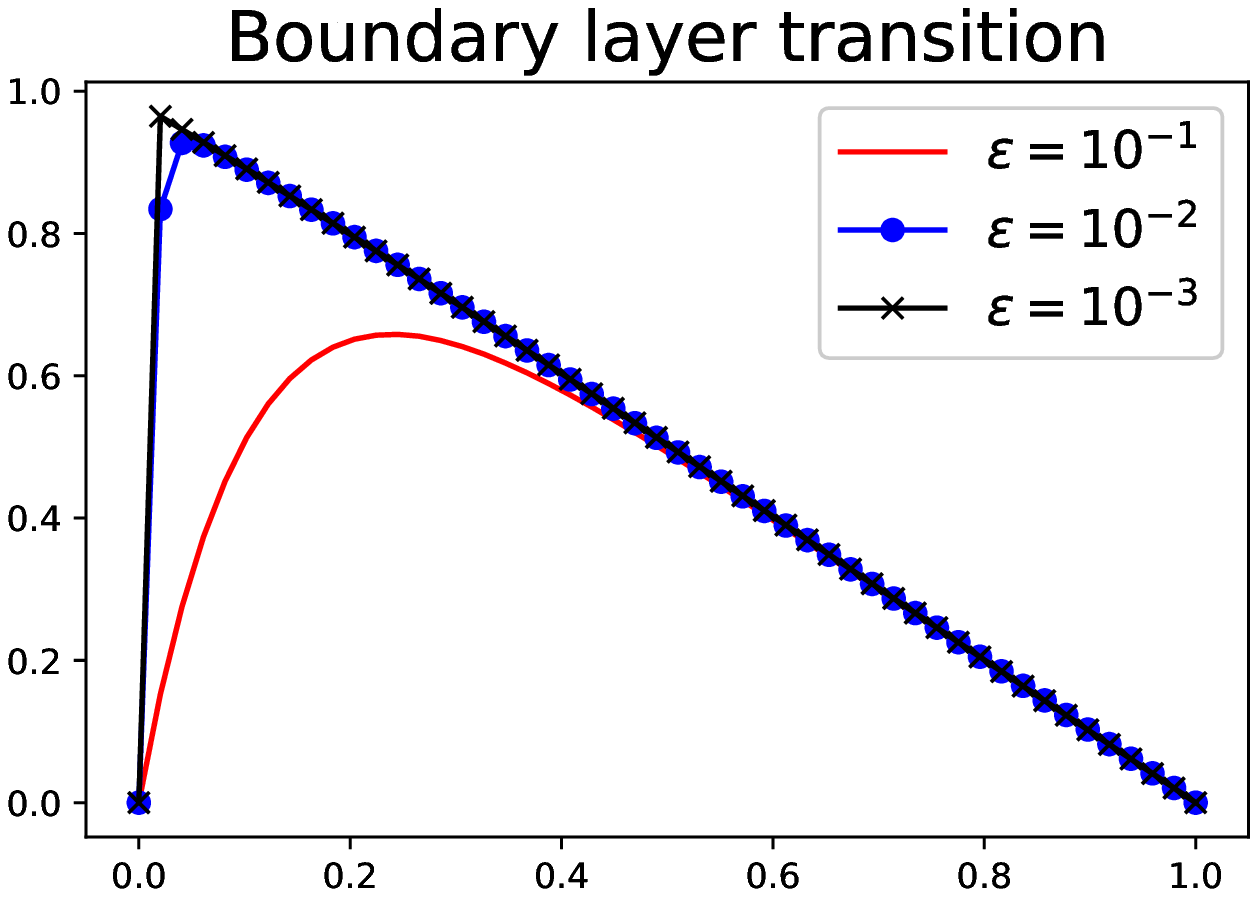}
    \caption{Boundary layer transition for \eqref{e:con_dif_eqn} with respect to different viscosity parameters, $\ep$, is displayed. As $\ep$ decreases, thickness of the boundary layer gets smaller.}
    \label{fig:2.3}
\end{figure}

Now, we compare the 
performance of 
the usual 2-layer (modified) PINN approximation $\overline{u}$ 
in (\ref{e:app_sol_mod}) 
and that of 
our new semi-analytic 2-layer PINN approximation $\widetilde{u}$:

The Fig.\ref{fig:classical_NN} 
shows that the usual PINN, constructed in Section \ref{S.reg} fails to 
approximate the solution of 
the  singularly perturbation (\ref{e:con_dif_eqn_limit}). 
The numerical results below 
in Fig. \ref{fig:2.1} 
and Table \ref{t:MAIN} 
confirm that 
the semi-analytic enriched 2-layer PINN  performs  much better than the usual PINN method, 
thanks to the corrector function embedded in the scheme. 
Note that 
our semi-analytic PINN enriched with corrector produces stable and accurate approximate  solutions, independent of the small parameter $\ep$, as  shown in the Fig. \ref{fig:2.3}. 
\begin{figure}[H]
\centering
     \begin{subfigure}[b]{0.4\textwidth}
         \centering
         \includegraphics[width=\textwidth]{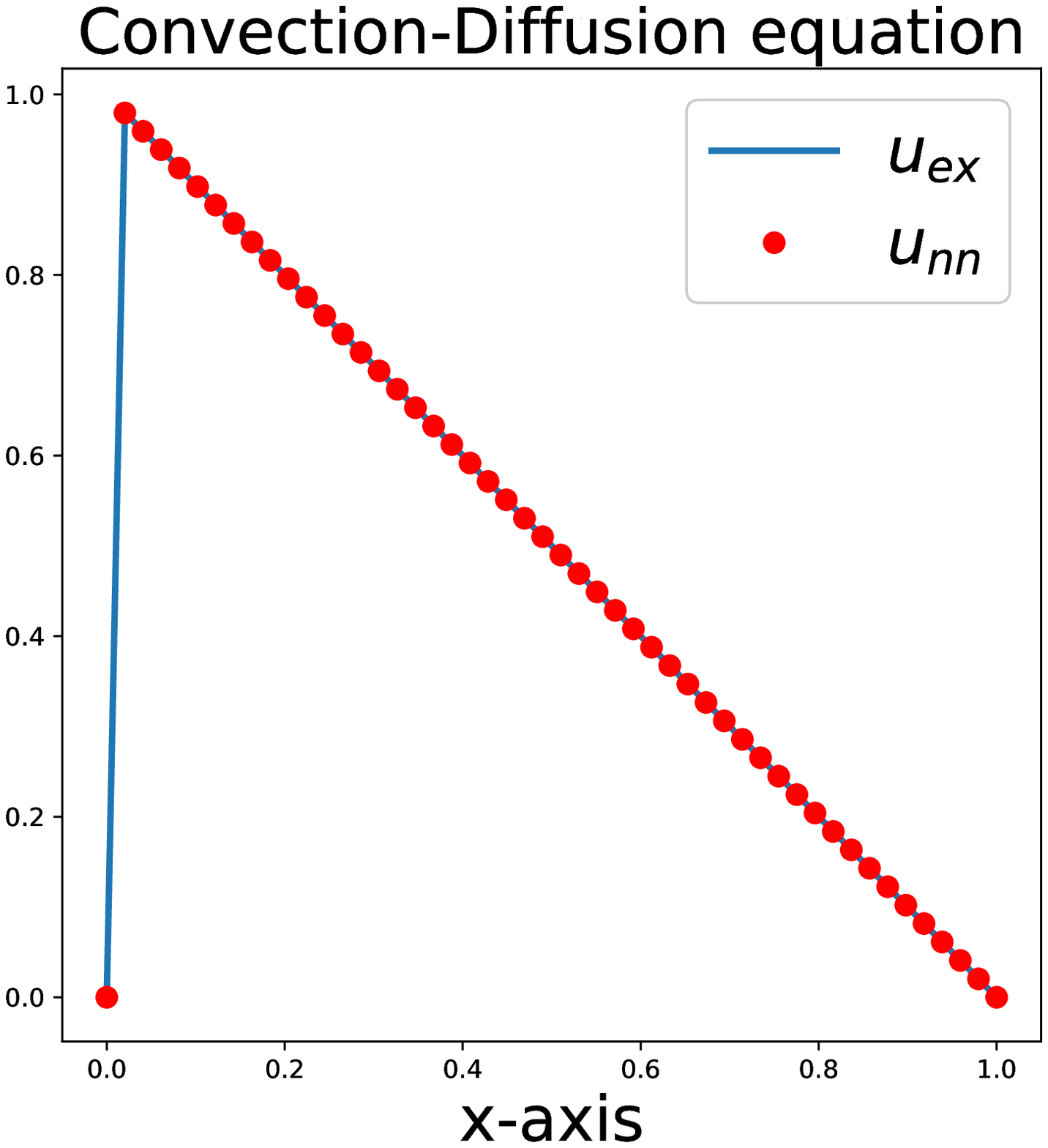}
         \caption{Exact sol. vs semi-analytic NN approx. for $f=1$}
     \end{subfigure}
\hspace{25mm}
     \begin{subfigure}[b]{0.4\textwidth}
         \centering
         \includegraphics[width=\textwidth]{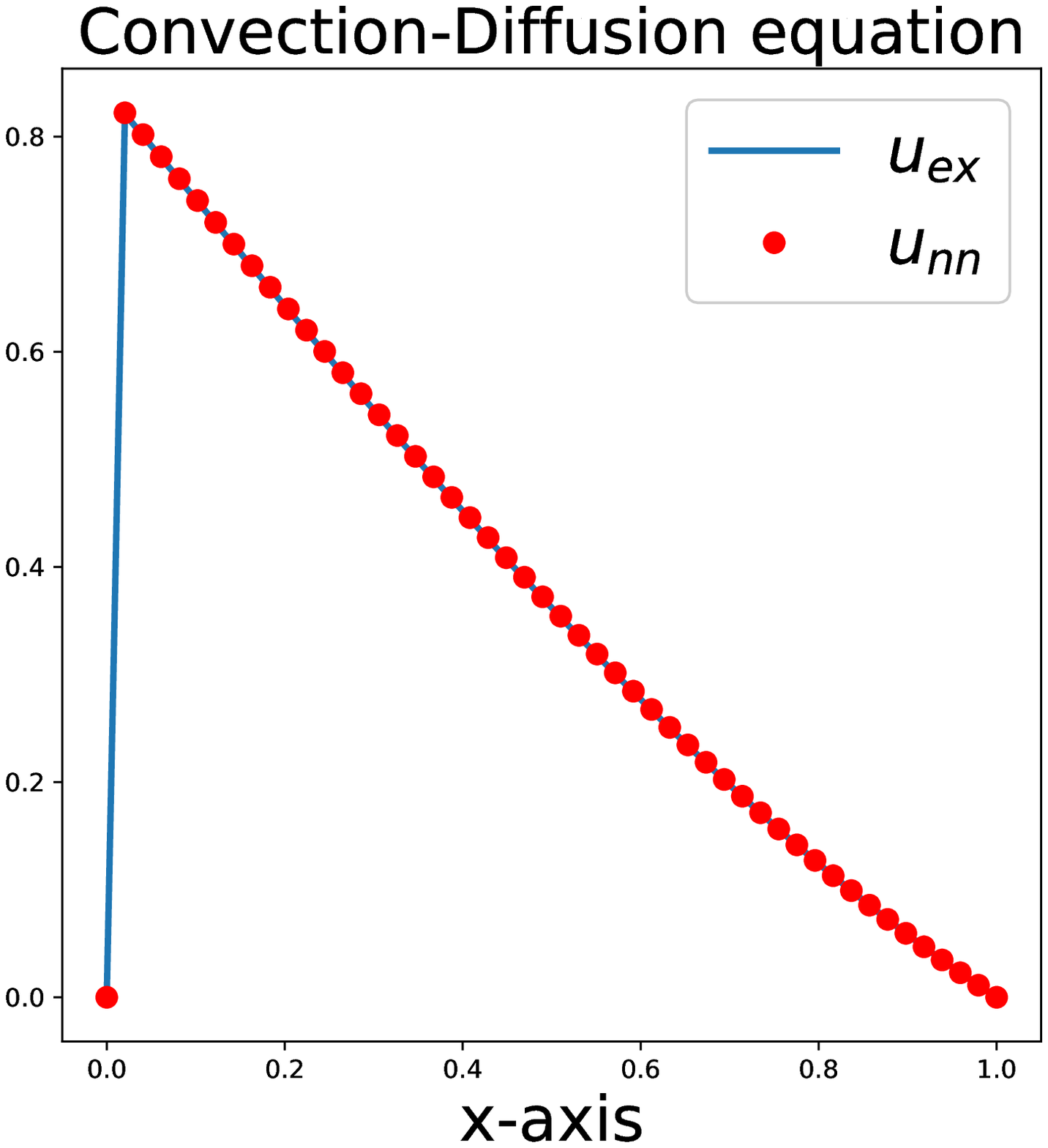}
         \caption{Exact sol. vs semi-analytic NN approx. for $f=\cos(x)$}
     \end{subfigure}
     \caption{The exact sol. and the approx. sol. of the convection-diffusion equation \eqref{e:con_dif_eqn} with $\epsilon = 10^{-4}$ are displayed. 
     Here, the semi-analytic enriched 2-layer PINN \eqref{e:app_sol_enriched_con-def} is used for the predicted sol.
     In panel (A), the predicted sol. is close to the corresponding exact sol. with relative $L^2$ error: $7.295 \times 10^{-3}$. 
     In panel (B), the predicted sol. is close to the corresponding exact sol. with relative $L^2$ error: $9.012 \times 10^{-3}$.}\label{fig:2.1}
\end{figure}

\subsection{Reaction-diffusion equations}\label{S.rea_dif_eqn}
We construct in this Section 
the semi-analytic PINN for the following singularly perturbed reaction-diffusion equation:
\begin{equation}\label{e:rea_dif_eqn}
    \left\{
    \begin{array}{rl}
        - \ep u^\ep_{xx} + u^\ep   = f(x), & 0 < x <1, \\
        u^\ep = 0, & x=0, 1.
    \end{array}
    \right.
\end{equation}

The corresponding limit, as $\ep \rightarrow 0$, $u^0$ 
of $u^{\ep}$ (solution of  (\ref{e:rea_dif_eqn})) is obtained by formally replacing $\ep$ by $0$ in (\ref{e:rea_dif_eqn}):
\begin{equation}\label{e:rea_dif_eqn_limit}
    u^0   = f.
\end{equation}

Performing the matching asymptotics,  we find that the size of boundary layer for (\ref{e:rea_dif_eqn}) is 
of order $\ep^{1/2}$, 
and it appears near both ends of the domain, i.e., 
near $x = 0$ and $x = 1$. 
Writing the asymptotic equation of the difference 
$u^\ep - u^0$, 
we find the equation for the corrector 
$\varphi \sim u^\ep - u^0$ as 
\begin{equation}\label{e:cor_eqn_rea-dif}
    \left\{
        \begin{array}{rl}
            \spacer
            -\ep \varphi_{xx} + \varphi  = 0, 
            &
            0 < x < 1, \\
            \varphi= - u^0 , 
            &
            x = 0, 1.
        \end{array}
    \right.
\end{equation}

By solving (\ref{e:cor_eqn_rea-dif}), we find that 
\begin{equation}\label{e:cor__rea-dif}
    \varphi(x)
        =
            -u^0(0) \, e^{-x/\sqrt{\ep}}
            -u^0(1) \, e^{-(1-x)/\sqrt{\ep}}
            + e.s.t., 
\end{equation}
i.e., 
the fast decaying part $\varphi$ of $u^\ep$ is 
the sum of 
two exponentially decaying functions from each part of the boundary points $x=0$ and $x=1$, 
scaled by the stretched variables $x/\sqrt{\ep}$ and 
$(1-x)/\sqrt{\ep}$ up to an exponentially small term.

Employing the simple energy estimates on the difference 
$u^\ep - (u^0 + \varphi)$, and 
then using the smallness of the corrector $\varphi$, 
we then 
further notice that 
\begin{equation}\label{e:con__rea-dif}
    \begin{array}{l}
    \spacer
        \|
            u^\ep - (u^0 + \varphi)
        \|_{L^2((0, 1))}
            \leq
                \kap \ep, \\
        \|
            u^\ep -  u^0  
        \|_{L^2((0, 1))}
            \leq
                \kap \ep^{\frac{1}{4}}, 
    \end{array}
\end{equation}
for a constant $\kap>0$ independent of $\ep$. 
Hence we notice here that 
the corrector $\varphi$ 
exhibits the singular behavior of $u^\ep$ at a small diffusivity $\ep$, i.e., 
the diffusive solution $u^\ep$ is decomposed into the sum of 
fast (decaying part) $\varphi$, 
and the slow part $u^0$.
\bigskip

Based on the boundary layer analysis above, 
using the corrector $\varphi$, 
we construct  
the new {\it semi-analytic enriched 2-layer PINN} for the problem (\ref{e:rea_dif_eqn}) in the form,
\begin{equation}\label{e:app_sol_enriched_rea-def}
    \widetilde{u}(x ; \, {\blds \theta}) 
        =
                \hat{u}({ {x}}; \, {\blds \theta})
            -
                \hat{u}({ {0}}; \, {\blds \theta}) \,
                e^{-x/\sqrt{\ep}}
            -
                \hat{u}({ {1}}; \, {\blds \theta}) \,
                e^{-(1-x)/\sqrt{\ep}}
            , 
\end{equation}
where $\hat{u}$ is exactly the same as in (\ref{e:NN_modified}). 
Because the effect of the exponentially decaying function $e^{-x/\sqrt{\ep}}$ 
(or $e^{-(1-x)/\sqrt{\ep}}$) 
on the boundary point at $x=1$ (or $x = 0$) 
is exponentially small with respect to the small $\ep$, 
the enriched PINN approximation $\widetilde{u}$ 
attains the zero boundary value at $x = 0,1$, 
up to an exponentially small 
(computationally negligible) 
error.

The loss for this problem (\ref{e:rea_dif_eqn}) is defined by 
\begin{equation}\label{e:loss_PINNs_mod_rea-dif}
    \mathcal{L}({\blds \theta}, \, \mathcal{T})
        =
            \frac{1}{|\mathcal{T} |}
            \sum_{ {\blds x}  \in \mathcal{T} }
            \|
            -\ep \widetilde{u}_{xx}
            - \widetilde{u} 
            - f 
            \|_{L^2(0, 1)}^2, 
\end{equation}
where 
the training set $\mathcal{T}$ is chosen as a set of scattered points in 
{\red $(0, 1)$}, and 
\begin{equation}\label{e:loss_eqn_rea-dif}
    -\ep \widetilde{u}_{xx}
    - \widetilde{u} 
    - f 
    =
            - \ep (x-1) \hat{u}_{xx} ({ {x}}; \, {\blds \theta})
            - \hat{u} ({ {x}}; \, {\blds \theta})
            - f
            , 
\end{equation}
with $\hat{u}_{xx}$ in (\ref{e:NN_modified_der}). 
In derivation of (\ref{e:loss_eqn_rea-dif}), 
we used the fact that 
the exponentially decaying functions 
$e^{-x/\sqrt{\ep}}$ 
and $e^{-(1-x)/\sqrt{\ep}}$  
satisfy the equation (\ref{e:rea_dif_eqn}).  
\bigskip
\bigskip

Now, we notice from Fig. \ref{fig:2.4} and Table \ref{t:MAIN}  
below that 
our semi-analytic 2-layer PINN, 
enriched with the corrector,  
approximates well the 
solution of (\ref{e:rea_dif_eqn}):

\begin{figure}[H]
    \centering
    \includegraphics[width=0.4\textwidth]{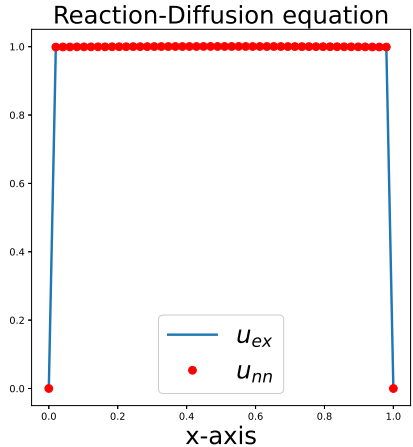}
    \caption{The exact sol. and the approx. sol. of the reaction-diffusion equation \eqref{e:rea_dif_eqn} with $f=1$ and $\epsilon = 10^{-4}$ are displayed. 
     The semi-analytic 2-layer PINN \eqref{e:app_sol_enriched_rea-def} is used for the predicted sol.
     The predicted sol. is close to the corresponding exact sol. with relative $L^2$ error: $1.245 \times 10^{-4}$.}
    \label{fig:2.4}
\end{figure}

\section{Enriched Neural Network: nonlinear singularly perturbed boundary value problems}\label{S.Sing_eqns_non}

In this section, we propose the 2-layer modified PINN for the nonlinear singularly perturbed boundary value problems. We first consider the convection-diffusion equations with
a non-linear reaction. The shape of the boundary layer profile is similar to that in the linear
case. We then examine the stationary Burgers’ equation with a small viscosity parameter. Since the analysis and computation of singularly perturbed Burgers’ equations are not
straightforward, our neural network requires careful numerical treatment.

\subsection{Convection-diffusion equations with a non-linear reaction}\label{S.rea_con_dif_eqn}
In this section, 
we apply our methodology of semi-analytic enriched PINNs to a certain non-linear equation 
for which 
we can determine the profile of boundary layer. 
Note that the corresponding boundary layer analysis is not  straightforward at all; 
see \cite{JP07}.

We consider the singularly perturbed convection-diffusion equation with a non-linear reaction term:
\begin{equation}\label{e:rea_con_dif_eqn_nonlinear}
    \left\{
    \begin{array}{rl}
        \spacer
        - \ep u^\ep_{xx} - u^\ep_x + (u^\ep)^3   = f(x), & 0 < x <1, \\
        u^\ep = 0, & x=0, 1.
    \end{array}
    \right.
\end{equation}
The corresponding limit problem at $\ep = 0$ is given by 
\begin{equation}\label{e:rea_con_dif_eqn_nonlinear_limit}
    \left\{
    \begin{array}{rl}
        \spacer
        - u^0_x + (u^0)^3   = f(x), & 0 < x <1, \\
        u^0 = 0, & x= 1.
    \end{array}
    \right.
\end{equation}
The well-posedness and the regularity of (\ref{e:rea_con_dif_eqn_nonlinear}) and (\ref{e:rea_con_dif_eqn_nonlinear_limit}) 
are well-studied and here we omit any further discussion on those issues; see, e.g., \cite{JP07} for the detailed information. 

Although it is a non-linear problem, 
the boundary layer associated with (\ref{e:rea_con_dif_eqn_nonlinear}) 
is linear. 
In fact, by performing the matching asymptotics, 
one can verify that 
the boundary layer of size $\ep$ occurs near the {\it outflow} boundary at $x = 0$, just like the linear convection-diffusion problem (\ref{e:con_dif_eqn}).  
Moreover we   
find that  
the asymptotic equation 
(with respect to the small $\ep$) for the corrector 
$\varphi$, 
which approximate the difference $u^\ep - u^0$, 
is given 
in the form,
\begin{equation}\label{e:cor_eqn_con-rea-dif}
    \left\{
        \begin{array}{rl}
            \spacer
            -\ep \varphi_{xx} - \varphi_x = 0, 
            &
            0 < x < 1, \\
            \varphi= - u^0 , 
            &
            x = 0, 
        \end{array}
    \right.
\end{equation}
and hence 
the corrector $\varphi$ is 
explicitly written in the form,
\begin{equation}\label{e:cor__con-rea-dif}
    \varphi(x)
        =
            -u^0(0) \, e^{-x/\ep}
            + e.s.t. 
\end{equation}
Concerning the detailed boundary layer analysis as well as the convergence results of 
$u^\ep$ to $u^0$, see, e.g., \cite{JP07}. 
\bigskip

By enriching the 2-layer PINNs with the profile of the corrector $\varphi$ above,  
we define 
the new {\it semi-analytic 2-layer PINN} for the problem (\ref{e:rea_con_dif_eqn_nonlinear}) in the form,
\begin{equation}\label{e:app_sol_enriched_con-rea-def}
    \widetilde{u}(x ; \, {\blds \theta}) 
        =
            (x-1) \, 
            ( 
                \hat{u}({ {x}}; \, {\blds \theta})
            -
                \hat{u}({ {0}}; \, {\blds \theta}) \,
                e^{-x/\ep}
            ), 
\end{equation}
with $\hat{u}$ as in (\ref{e:NN_modified}).

The loss for this problem is defined by 
\begin{equation}\label{e:loss_PINNs_mod_con-rea-dif}
    \mathcal{L}({\blds \theta}, \, \mathcal{T})
        =
            \frac{1}{|\mathcal{T} |}
            \sum_{ {\blds x}  \in \mathcal{T} }
            \|
            -\ep \widetilde{u}_{xx}
            - \widetilde{u}_x
            + \widetilde{u}^3
            - f 
            \|_{L^2(0, 1)}^2, 
\end{equation}
where 
the training set $\mathcal{T}$ is chosen as a set of scattered points in 
{$(0, 1)$}. 
The derivatives of our enriched 2 layer approximation $\widetilde{u}$ are exactly the same as in (\ref{e:2layer_app_enriched_der}), 
and thus, using the fact that 
the exponentially decaying function $e^{-x/\ep}$ 
satisfies the equation (\ref{e:cor_eqn_con-rea-dif}), 
we find for 
the loss function  
(\ref{e:loss_PINNs_mod_con-rea-dif}) that  
\begin{equation}\label{e:loss_eqn_con-rea-dif}
\begin{array}{l}
\spacer
    -\ep \widetilde{u}_{xx}
    - \widetilde{u}_x
    + \widetilde{u}^3
    - f \\
    \spacer
    \qquad
        =
            - \ep (x-1) \hat{u}_{xx} ({ {x}}; \, {\blds \theta})
            - (x - 1 + 2\ep) \hat{u}_{x} ({ {x}}; \, {\blds \theta})
            - \hat{u} ({ {x}}; \, {\blds \theta})
            -
                \hat{u} (0; \, {\blds \theta}) e^{-x/\ep}
                \\
    \qquad \quad
            \,
            +
            (x-1)^3 \, 
            ( 
                \hat{u}({ {x}}; \, {\blds \theta})
            -
                \hat{u}({ {0}}; \, {\blds \theta}) \,
                e^{-x/\ep}
            )^3
            - f
            .
\end{array}
\end{equation}
Here the derivatives of $\hat{u}$ are given in (\ref{e:NN_modified_der}).  
Note that 
all the terms in (\ref{e:loss_eqn_con-rea-dif}) stay bounded as 
$\ep$ vanishes, 
and hence 
our semi-analytic enriched 2 layer PINN 
produces 
an 
accurate 
approximation $\widetilde{u}$ for (\ref{e:rea_con_dif_eqn_nonlinear}),  
independent of the small parameter $\ep$.
\bigskip

We observe from Fig. \ref{fig:2.5} and Table \ref{t:MAIN}  
below that 
our semi-analytic 2-layer PINN, 
enriched with the corrector,  approximates well the 
solution of (\ref{e:rea_con_dif_eqn_nonlinear}):

\begin{figure}[H]
    \centering
    \includegraphics[width=0.4\textwidth]{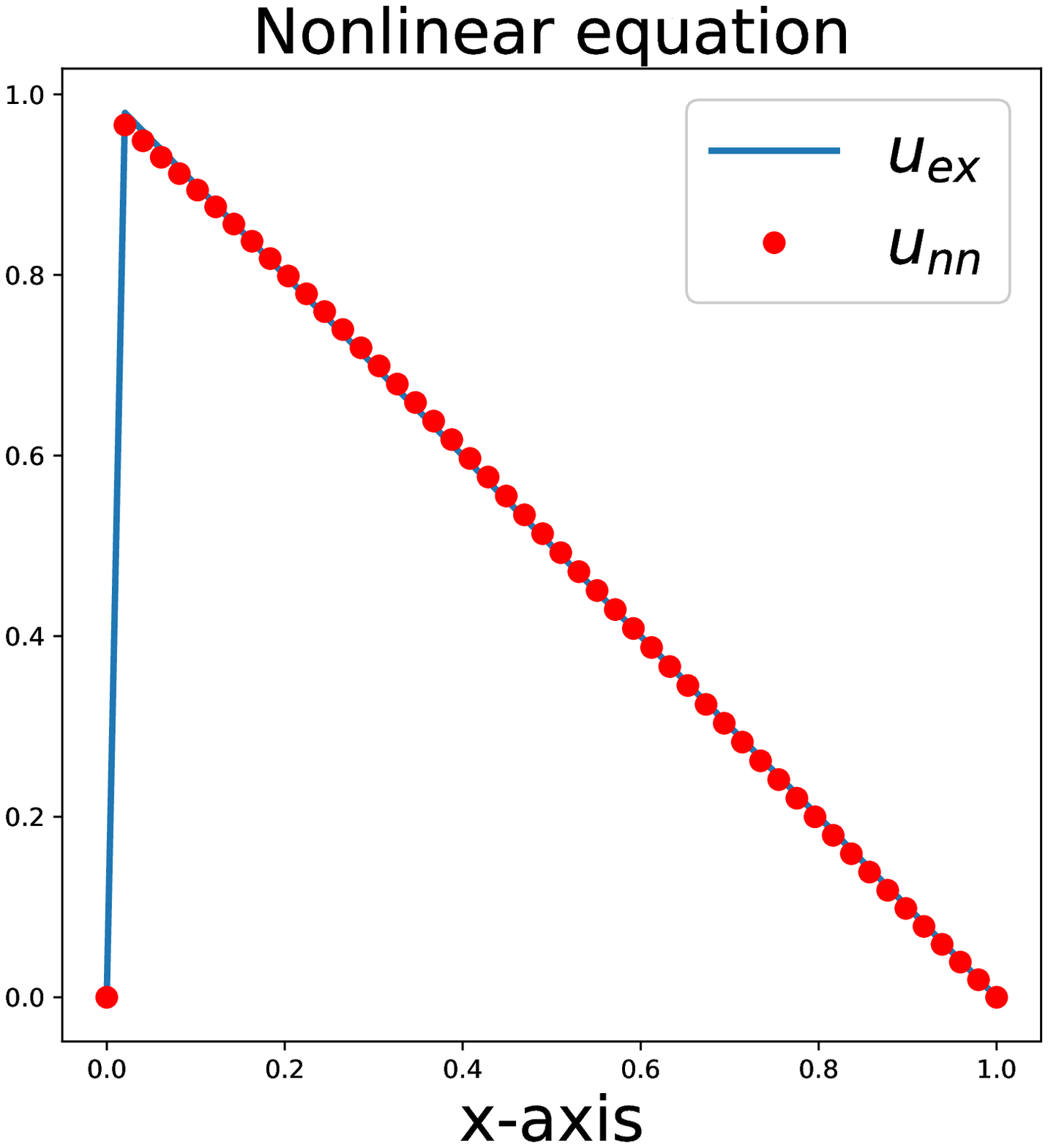}
    \caption{The exact sol. and the approx. sol. of \eqref{e:rea_con_dif_eqn_nonlinear} with $f=1$ and $\epsilon = 10^{-4}$ are displayed. 
     Here the semi-analytic 2-layer PINN \eqref{e:app_sol_enriched_con-rea-def} is used for the predicted sol.
     The predicted sol. is close to the corresponding exact sol. with relative $L^2$ error: $5.474 \times 10^{-2}$.}
    \label{fig:2.5}
\end{figure}

\subsection{Stationary Burgers' equation}\label{S.Burgers_eqn}
In this section, 
we consider the stationary 1D Burgers' equation in a bounded interval $(0, 1)$ as  
\begin{equation}\label{e:Burgers}
\left\{\begin{array}{rl}
	\spacer
	-\ep u^\ep_{xx} + u^\ep u^\ep_x & = f(x), 
		 \quad 0 < x < 1,\\
        u^\ep & = -\alpha, \quad  x=0,\\
        u^\ep & = -\beta,  \quad x=1.
\end{array}\right.
\end{equation}
Here  
$0 < \epsilon \ll 1$ is a small viscosity parameter, 
$f$ is a smooth data, independent of $\ep$, 
and $\alpha$, $\beta$ are positive constants. 
We set the boundary values of $u^\ep$ at $x=0, 1$ 
as negative numbers so that 
$u^\ep < 0$ for all $0 \leq x \leq 1$.   
Hence, consequently,  
the convection $u^\ep u^\ep_x $ occurs always in one direction from right to left 
and the boundary layer occurs near only   
the out-flow boundary at $x=0$.

For the sake of convenience on our computations below, we set $\alpha = 1$ and $\beta = 1$. 
We also assume that the data $f$ satisfies the following condition so that 
the  limit solution $u^0$ is well-defined and explicitly written as in (\ref{e:u^0}) below:
\begin{equation}\label{e:f_condition_Burgers}
    \int_1^x f(x^\pri) \, dx^\pri
        > - \dfrac{1}{2}, 
    \quad 
    \text{for any } 0 \leq x \leq 1; 
\end{equation}
see, e.g., \cite{CJL19, GJL22} for more information and full boundary layer analysis for this version of stationary 1D Burgers' equation.  
 
The corresponding limit (inviscid) problem is obtained 
by setting $\ep = 0$ in (\ref{e:Burgers}) 
and imposing the in-flow boundary condition at $x=1$:
\begin{equation}\label{e:Burgers_limit}
	\left\{
	\begin{array}{l}
		\spacer
			u^0 u^0_x = f(x), \quad 0 < x < 1, \\
			u^0(1) = - 1.
	\end{array}
	\right.
\end{equation}
The formal limit $u^0$, a solution of (\ref{e:Burgers_limit}), is given in the form, 
\begin{equation}\label{e:u^0}
    u^0 (x) 
        =
            - \Big(
                2 \int_1^x f(x^\pri) \, dx^\pri 
                + 1
            \Big)^{\frac{1}{2}}, 
\end{equation}
and hence we infer from (\ref{e:f_condition_Burgers}) that 
\begin{equation}\label{e:u^0_sign}
    u^0 < 0 \quad 
    \text{for any } 0 \leq x \leq 1.
\end{equation} 

Performing the matching asymptotics, 
we find that the size of boundary layer for (\ref{e:Burgers}) is 
of order $\ep$, 
and it appears near the out-flow boundary at $x = 0$.  
Writing the asymptotic equation of the difference 
$u^\ep - u^0$, 
we find the following ({\it non-linear}) 
asymptotic equation for   
$\varphi \sim u^\ep - u^0$ as 
\begin{equation}\label{e:cor_eqn_Burgers}
    \left\{
        \begin{array}{rl}
            \spacer
            -\ep \varphi_{xx} 
            + 
            \big(
            u^0(0) + \varphi 
            \big)
            \varphi_x  = 0, 
            &
            0 < x < 1, \\
            \spacer
            \varphi= -1 - u^0 , 
            &
            x = 0, \\
            \varphi \rightarrow 0
            &
            \text{as }
            x \rightarrow 1.
        \end{array}
    \right.
\end{equation}
By integrating (\ref{e:cor_eqn_Burgers}) from $x^\pri = 1$ to $x^\pri = x$, 
we write the first order equation,  
\begin{equation}\label{e:cor_eqn_Burgers_1st}
    \left\{
        \begin{array}{rl}
            \spacer
            -\ep \varphi_{x} 
            + 
            u^0(0) \varphi 
            +
            \dfrac{1}{2}
            \varphi^2  = 0, 
            &
            0 < x < 1, \\
            \varphi= -1 - u^0 , 
            &
            x = 0.
        \end{array}
    \right.
\end{equation}
Then, solving the equation above, we find  
\begin{equation}\label{e:cor_Burgers}
    \varphi(x)
        =
            \dfrac{
                2u^0(0) \big(1 + u^0(0)\big)
                e^{u^0(0)x/{\ep}}
            }
            {
                1 - u^0(0) 
                -
                \big(1 + u^0(0)\big)
                e^{u^0(0)x/{\ep}}
            }
            .
\end{equation}
The fast decaying part of $u^\ep$ (near the out-flow boundary at $x=0$) is 
hence described by the corrector $\varphi$ above. 
Moreover, by performing the boundary layer analysis as in, e.g., \cite{CJL19, GJL22}, one can verify that 
\begin{equation}\label{e:conv_bur}
    \|
        u^\ep - (u^0 + \varphi )
    \|_{L^2(0, 1)}
        \leq
            \kap \ep,
\end{equation}
for a constant $\kap >0$ independent of $\ep$, and hence 
obtain the vanishing viscosity limit as well:
\begin{equation}\label{e:VVL_bur}
    \|
        u^\ep - u^0  
    \|_{L^2(0, 1)}
        \leq
            \kap \ep^\frac{1}{2}.
\end{equation}

\bigskip

Now, thanks to the asymptotic analysis above, 
we construct below the semi-analytic 2-layer PINNs enriched by
the profile of the corrector $\varphi$:

We first normalize the boundary value of $\varphi$ in (\ref{e:cor_Burgers}), 
and introduce 
the normalized corrector $\widetilde{\varphi}$, 
which describe the boundary layer profile for (\ref{e:Burgers}), 
in the form, 
\begin{equation}\label{e:cor_Burgers_nor}
    \widetilde{\varphi}(x)
        =
            \dfrac{{\varphi}(x)}{{\varphi}(0)}
            .
\end{equation}
Then, we define 
our new {\it semi-analytic enriched} 2-layer PINN for the problem (\ref{e:Burgers}) as 
\begin{equation}\label{e:app_sol_enriched_Burgers}
\begin{split}
    \widetilde{u}(x ; \, {\blds \theta}) 
        & =
            (x-1) \,  \hat{u}({ {x}}; \, {\blds \theta})
             + \widetilde{\varphi}(x) \hat{u}({ {0}}; \, {\blds \theta}) - g(x)\\
        & = \Big( \text{by setting $\widetilde{\varphi}^*(x) := \widetilde{\varphi}(x) \hat{u}({ {0}}; \, {\blds \theta})$} \Big)\\
        & =   (x-1) \,  \hat{u}({ {x}}; \, {\blds \theta})
             + \widetilde{\varphi}^*(x) - g(x),
\end{split}             
\end{equation}
where $\hat{u}$ is exactly the same as in (\ref{e:NN_modified}), which 
was used for the usual 2-layer PINNs approximation $\overline{u}$ 
in (\ref{e:app_sol_mod}), and $g(x)$ is a simple boundary lifting function,
\begin{equation}
    g(x) = (\beta - \alpha) x +\alpha = 1, 
    \quad 
    (\text{by using $\alpha = \beta = 1$});
\end{equation}
note that the choice of this lifting $g = 1$ is for our convenience, but any other lifting, which gives the value $1$ at $x=0,1$, produces the same computational results as those we obtain in this article.

The loss for this problem is defined by 
\begin{equation}\label{e:loss_PINNs_mod_bur}
    \mathcal{L}({\blds \theta}, \, \mathcal{T})
        =
            \frac{1}{|\mathcal{T} |}
            \sum_{ {\blds x}  \in \mathcal{T} }
            \|
            -\ep \widetilde{u}_{xx} + \widetilde{u} \widetilde{u}_x 
            - f
            \|_{L^2(0, 1)}^2, 
\end{equation}
where 
the training set $\mathcal{T}$ is chosen as a set of scattered points in 
{$(0, 1)$}. 
The derivatives of our enriched 2 layer approximation $\widetilde{u}$ are given by 
\begin{equation}\label{e:2layer_app_enriched_der_bur}
\begin{array}{l}
    \spacer
    \dfrac{d}{dx} \widetilde{u}
        =
            \hat{u}({ {x}}; \, {\blds \theta})
            + (x-1) \hat{u}_x({ {x}}; \, {\blds \theta})
             + \widetilde{\varphi}^*_x(x)\\
    \dfrac{d^2}{dx^2} \widetilde{u}
        =
            2 
                \hat{u}_x({ {x}}; \, {\blds \theta})
            +
            (x-1)
                \hat{u}_{xx}({ {x}}; \, {\blds \theta})
                + \widetilde{\varphi}^*_{xx}(x).
\end{array}
\end{equation}
{  
We recall from (\ref{e:conv_bur}) that 
the asymptotic expansion of $u^\ep$ is 
well defined in the sense that 
\begin{equation}\label{e:asymp_exp_bur}
u^\ep 
    = 
        u^0 + \varphi + \mathcal{O}(\ep).
\end{equation}
Then, because the corrector $\varphi$ satisfies  
the equation  (\ref{e:cor_eqn_Burgers})$_1$, and because we use $\widetilde{u}$ to approximate $u^\ep$, i.e., 
$u^\ep \simeq \widetilde{u}$, 
we observe that the corrector $\varphi$ 
satisfies 
\begin{equation}\label{eqn_cor_bur_app_TEMP}
    -\ep \varphi_{xx} 
            + 
            \big(
            \widetilde{u}
            -
            (u^0 - u^0(0)
            )
            \big)
            \varphi_x  = 
            \mathcal{O}(\ep), 
\end{equation}
which is equivalent, 
in terms of the normalized corrector, 
to 
\begin{equation}\label{eqn_cor_bur_app_TEMP2}
    -\ep \widetilde{\varphi}^*_{xx} 
            + 
            \widetilde{u} \,
            \widetilde{\varphi}^*_x
            -
            (u^0 - u^0(0))
            \widetilde{\varphi}^*_x  
            = 
            \mathcal{O}(\ep). 
\end{equation}
Finally, 
using (\ref{e:2layer_app_enriched_der_bur}) and 
the fact that 
the normalized corrector 
$\widetilde{\varphi}$ 
satisfies the equation  (\ref{eqn_cor_bur_app_TEMP2}) above, 
we find for 
the loss function  
(\ref{e:loss_PINNs_mod_bur}) that  
\begin{equation}\label{e:loss_eqn_bur}
\begin{array}{rl}
    -\ep \widetilde{u}_{xx} + \widetilde{u} \widetilde{u}_x 
            - f
            =
            &
            \spacer
                - 2\ep
                \hat{u}_x
                - \ep (x-1) \hat{u}_{xx}
                + 
                \big(
                    (x-1)  \hat{u} 
                \big)_x
                \big(
                    (x-1) \hat{u} - \widetilde{\varphi}^* - g
                \big)\\
            &    +
                \big(
                    u^0 - u^0(0)
                \big)\widetilde{\varphi}^*_x ,
\end{array}
\end{equation}
where the derivatives of $\hat{u}$ are given in (\ref{e:NN_modified_der}).  
Because  
$u^0 - u^0(0) \simeq \mathcal{O}(\ep)$ 
and 
$\widetilde{\varphi}_x \simeq \ep^{-1} e^{-x/\ep}$ from (\ref{e:cor_Burgers}), 
we notice that 
all the terms in (\ref{e:loss_eqn_bur}) stay bounded as 
$\ep$ vanishes, 
and hence 
our semi analytic 2 layer PINN 
produces 
an 
accurate 
approximation $\widetilde{u}$ for (\ref{e:Burgers}),  
independent of the small parameter $\ep$.
}

To measure performance of the semi-analytic 2-layer PINN approximation, we manufacture an exact solution to the Burgers' equations. 
For this, we use a numerical solution with a large number of discretization point since the exact solution is not available in general.
More precisely, we implement the Burgers' equations, \eqref{e:Burgers}, with $f = -1$ using the spectral element method with the number of elements $M = 2048$.
Since the collocation points (input points) of the neural network is much less than $M=2048$, we exploit spline interpolation when comparing the manufactured solution with the predicted solution.
Since a lot of collocation points are used in the spline interpolation, the numerical error from the interpolation is much less than the approximation error from the semi-analytic enriched 2-layer PINN approximation.

\begin{figure}[H]
    \centering
    \includegraphics[width=0.4\textwidth]{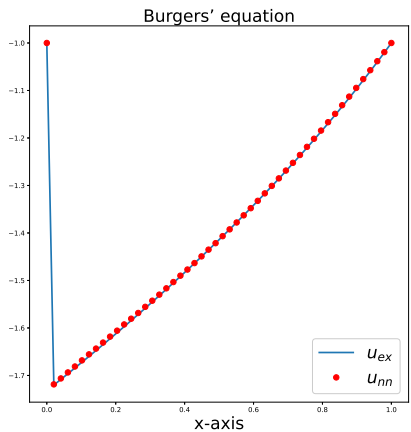}
    \caption{The exact sol. and the approx. sol. of \eqref{e:Burgers} with $f=-1$ and $\epsilon = 10^{-4}$ are displayed. 
     Here the semi analytic 2-layer PINN \eqref{e:app_sol_enriched_Burgers} is used for the predicted sol.
     The predicted sol. is close to the corresponding exact sol. with relative $L^2$ error: $2.599 \times 10^{-3}$.}
    \label{fig:2.6}
\end{figure}

\begin{figure}[H]
    \centering
    \includegraphics[width=0.5\textwidth]{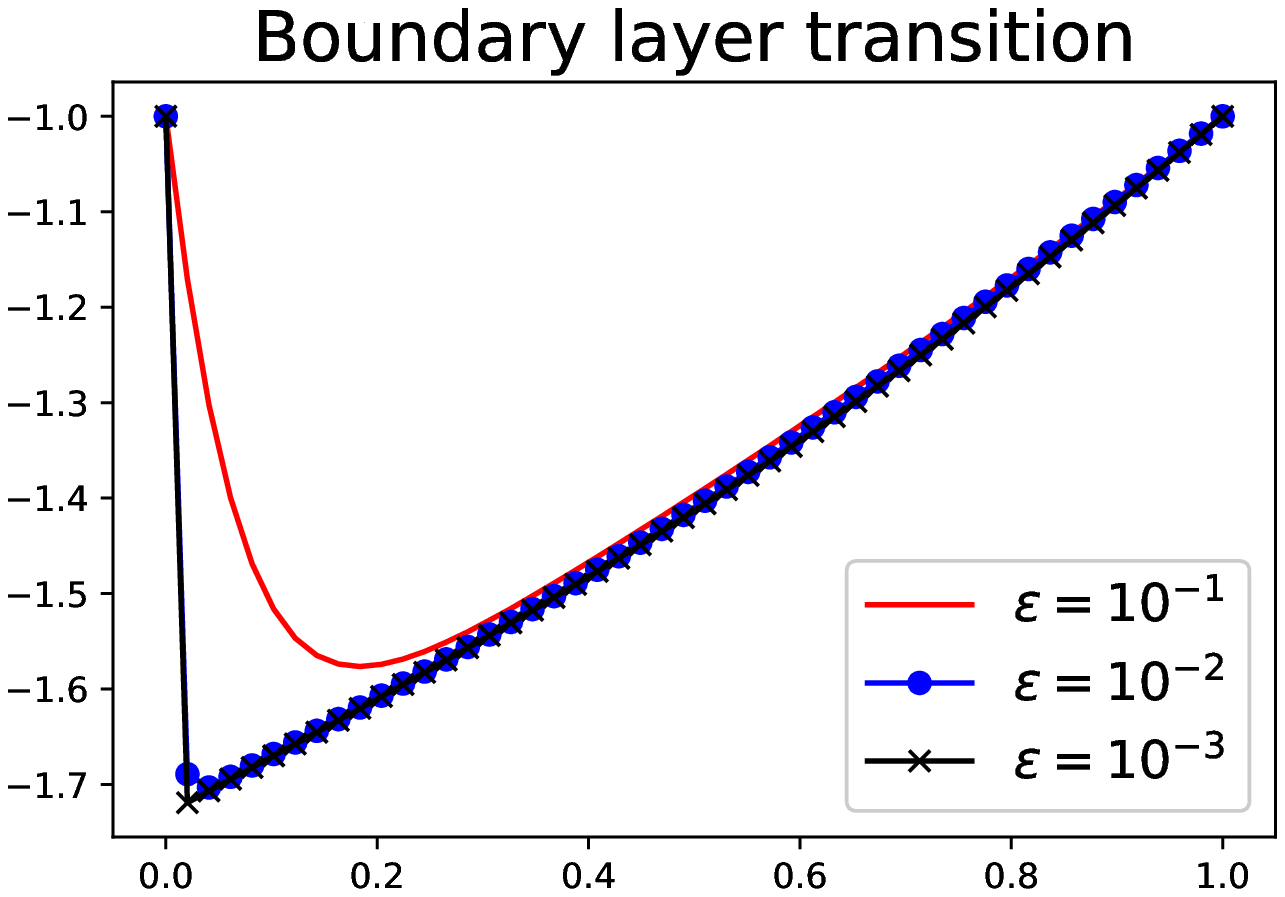}
    \caption{Boundary layer transition for \eqref{e:Burgers} with respect to different viscosity parameters, $\ep$, is displayed. As $\ep$ decreases, thickness of the boundary layer becomes smaller.}
    \label{fig:2.7}
\end{figure}

We observe from Fig. \ref{fig:2.6} and \ref{fig:2.7}, 
and Table \ref{t:MAIN}  
below that 
our semi-analytic 2-layer PINN, 
enriched with the corrector, 
approximates well the 
solution of (\ref{e:Burgers}). 
Our semi-analytic PINN   produces stable and accurate approximate solutions, independent of the small parameter $\ep$, as  shown in 
the Fig. \ref{fig:2.7}.

\begin{table}[H]
{\small
\begin{tabular}{|l|l|l|l|l|l|}
\hline
 N
 &  ECD ($\ep = 10^{-4}$) & CCD ($\ep = 10^{-4}$) & LRD ($\ep = 10^{-8}$) & NCD ($\ep = 10^{-4}$) & BE ($\ep = 10^{-4}$) \\ \hline
 $50$ & $7.295 \times 10^{-3}$ & $9.977 \times 10^{-1}$ & $1.245 \times 10^{-4}$ & $5.474 \times 10^{-2}$ & $2.599 \times 10^{-3}$ \\ \hline
 $100$ & $3.724 \times 10^{-3}$ & $9.931 \times 10^{-1}$  & $5.382 \times 10^{-4}$ & $ 4.885 \times 10^{-3}$ & $1.982 \times 10^{-3}$  \\ \hline
 $200$ & $2.148 \times 10^{-3}$  & $9.813 \times 10^{-1}$  & $4.245 \times 10^{-4}$ &  $ 5.735 \times 10^{-3}$ & $1.474 \times 10^{-3}$\\ \hline
 $400$ & $9.310 \times 10^{-4}$ & $9.550 \times 10^{-1}$ & $5.317 \times 10^{-4}$ & $8.912 \times 10^{-4}$ & $9.422 \times 10^{-4}$ \\ \hline
\end{tabular} \vspace{1mm}
}
\caption{Relative $L^2$ errors. Here, CD, LRD, NCD, and BE stand for, in turns, the convection-diffusion, linear reaction-diffusion, nonlinear convection-diffusion equations, and Burgers' equation  respectively. 
For the CD equation, 
ECD and CCD means respectively 
the semi-analytic enriched PINN and 
the 2 layer PINN without using the corrector. 
We set $N$ to be the number of data points. In this table, the maximum number of iterations is $50,000$, but early stopping is used to avoid over-fitting.}
\label{t:MAIN}
\end{table}
 
\section{Conclusion}\label{S.Con}
In this work, we have presented 
a semi-analytic approach to 
improve the 
numerical performance of 
the 2-layer PINNs,    
applied to various 
singularly perturbed boundary value problems. 
For each singular perturbation problem under consideration, in particular,  including the non-linear Burgers' equation, 
we first derived the so-called {\it corrector} function, 
which is an analytic approximation  of the fast (stiff) part of the solution to each example inside the boundary layer. 
By embedding the correctors into the structure of  2-layer PINNs, 
we resolve the stiffness nature of approximate solutions  and 
build 
our new semi-analytic PINNs enriched by the correctors. 
Performing numerical simulations, 
we verify that 
our new semi-analytic enriched PINNs 
produce stable and convergent 
approximations of the solutions  
to all the singular perturbations considered in this article.

\section*{Acknowledgments}

\noindent 
Gie was partially supported by 
Ascending Star Fellowship, Office of EVPRI, University of Louisville; 
Simons Foundation Collaboration Grant for Mathematicians; 
Research R-II Grant, Office of EVPRI, University of Louisville; 
Brain Pool Program through the National Research Foundation of Korea (NRF) (2020H1D3A2A01110658).  
The work of Y. Hong was supported by Basic Science Research Program through the National Research Foundation of Korea (NRF) funded by the Ministry of Education (NRF-2021R1A2C1093579) and the Korea government(MSIT)(No. 2022R1A4A3033571). 
Jung was supported by the Basic Science Research Program through the National Research Foundation of Korea funded by the Ministry of Education (2018R1D1A1B07048325)


\begin{thebibliography}{100}


 


 
\bibitem{BN18}
J. Berg and K. Nystrom. 
\newblock {A unified deep artificial neural network approach to partial differential equations in complex geometries}.  
{\em Neurocomputing}, 
317 (2018), pp. 28 --41.

\bibitem{BE21}
J. Blechschmidt, O. G. Ernst.  
Three ways to solve partial differential equations with neural networks--a review.  
{\em GAMM-Mitteilungen}, 44 (2) (2021). 
 
 



\bibitem{CF21}
Yuan Cao, Zhiying Fang, Yue Wu, Ding-Xuan Zhou, and Quanquan Gu.
Towards Understanding the Spectral Bias of Deep Learning.
{\em Proceedings of the Thirtieth International Joint Conference on Artificial Intelligence}, 8 (2021). 



 
 
\bibitem{CX21}
Z. Chen, D. Xiu.  
On generalized residual network for deep learning of unknown dynamical systems.  
{\em Journal of Computational Physics}, 438 (2021). 


\bibitem{CCWX22}
Z. Chen, V. Churchill, K. Wu, D. Xiu.  
Deep neural network modeling of unknown partial differential equations in nodal space.   
{\em Journal of Computational Physics} 449 (2022).



\bibitem{CJL19}
	Junho Choi, Chang-Yeol Jung, and Hoyeon Lee.
	\newblock On boundary layers for the {B}urgers equations in a bounded
              domain.
	\newblock {\em Commun. Nonlinear Sci. Numer. Simul.}, 67:637--657, 2019.
	
 
\bibitem{EY18}
W. E and B. Yu. 
\newblock {The deep Ritz method: a deep learning-based numerical algorithm for solving variational problems}. 
{\em Communications in Mathematics and Statistics}, 
6 (2018), pp. 1--12.

\bibitem{GJL1} 
	\textnormal{G.-M. Gie, C.-Y. Jung, and H. Lee}.  
	{Enriched Finite Volume approximations of the plane-parallel flow at a small viscosity.}  
	\textit{Journal of Scientific Computing}, 84, 7 (2020).
	
\bibitem{GJL2} 
	{\textnormal{G.-M. Gie, C.-Y. Jung, and H. Lee}.  
	{Semi-analytic time differencing methods for singularly perturbed initial value problems.}
	{\em 
	Numerical Methods for Partial Differential Equations,} 
	38, 5, 
	1367 - 139, 2022.}
	
	
\bibitem{GJL22} 
	{\textnormal{G.-M. Gie, C.-Y. Jung, and H. Lee}.  
	{Semi-analytic shooting methods for Burgers' equation.} 
	{\em 
	Accepted in 
	Journal of Computational and Applied Mathematics.
	}
	}
	

\bibitem{Book16} 
	G.-M. Gie, M. Hamouda, C.-Y. Jung, and R. Temam, 
	\newblock {\em Singular perturbations and boundary layers}, volume~200 of {\em Applied Mathematical Sciences}.
	\newblock  Springer Nature Switzerland AG, 2018.  
	\newblock  https://doi.org/10.1007/978-3-030-00638-9






\bibitem{HK82} 
	\textnormal{H. Han and R. B. Kellogg}, 
	\textit{A method of enriched subspaces for the numerical solution of a parabolic singular perturbation problem.} 
	In: Computational and Asymptotic Methods for Boundary and Interior Layers, Dublin, pp.46-52 (1982).


\bibitem{HL21}
J. Han, Y. Lee.  
Hierarchical learning to solve partial differential equations using physics-informed neural networks.  
{\em arXiv preprint}, arXiv:2112.01254 (2021).



\bibitem{Ho95} 
	\textnormal{M. H. Holmes}, 
	\textit{Introduction to perturbation methods}, 
	Springer, New York, 1995.
	
\bibitem{HJL13}
	Youngjoon Hong, Chang-Yeol Jung, and Jacques Laminie.
	\newblock Singularly perturbed reaction-diffusion equations in a circle with numerical applications.
	\newblock {\em Int. J. Comput. Math.}, 90(11):2308--2325, 2013.


\bibitem{HJT14}
	Youngjoon Hong, Chang-Yeol Jung, and Roger Temam.
	\newblock On the numerical approximations of stiff convection-diffusion equations in a circle.
	\newblock {\em Numer. Math.}, 127(2):291--313, 2014.




\bibitem{JP07}
Chang-Yeol Jung and Du Pham. 
\newblock {{Singular perturbation of semi-linear reaction-convection equations in a channel and numerical applications}}.
\newblock {\textit{Advances in Differential Equations}}. 
\newblock {Vol.12, no. 3, 2007.} 


\bibitem{KKLPWY21}
G. E. Karniadakis, I. G. Kevrekidis, L. Lu, P. Perdikaris, S. Wang,
L. Yang.  
Physics-informed machine learning.  
{\em Nature Reviews Physics}, 3 (6) (2021) 422--440.




\bibitem{KZK19}
E. Kharazmi, Z. Zhang, G. E. Karniadakis.  
Variational physics-informed neural networks for solving partial differential equations.  
{\em arXiv preprint}, arXiv:1912.00873 (2019).


\bibitem{KZK21}
E. Kharazmi, Z. Zhang, and G. E. Karniadakis.  
\newblock {hp-VPINNs: Variational physics-informed neural networks with
domain decomposition}. 
{\em Comput. Methods in Appl. Mech. Eng.}, 
374 (2021).  

\bibitem{KDJH21}
S. Kollmannsberger, 
D. D’Angella, 
M. Jokeit, 
L. Herrmann.  
\newblock {
\it Deep Learning in Computational Mechanics}. 
{Studies in Computational Intelligence, Springer International Publishing, Cham (2021)}.  
https://doi.org/10.1007/978-3-030-76587-3


\bibitem{LDKRKS20}
L. Lu, M. Dao, P. Kumar, U. Ramamurty, G. E. Karniadakis, S. Suresh.  
Extraction of mechanical properties of materials through deep learning
from instrumented indentation.  
{\em Proceedings of the National Academy of Sciences}, 117 (13) (2020) 7052--7062.


\bibitem{LLF98}
I.E. Lagaris, A. Likas, and D.I. Fotiadis.
\newblock {Artificial neural networks for solving ordinary and partial differential equations}. 
{\em IEEE Trans Neural Netw.},  1998;9(5):987-1000. 

\bibitem{LMMK21}
L. Lu, X. Meng, Z. Mao, G. E. Karniadakis.  
Deepxde: A deep learning library for solving differential equations.  
{\em SIAM Review}, 63 (1) (2021) 208--228.


\bibitem{MZK20}
X. Meng, Z. Li, D. Zhang, G. E. Karniadakis.  
Ppinn: Parareal physics-informed neural network for time-dependent pdes.  
{\em Computer Methods in Applied Mechanics and Engineering}, 370 (2020).


\bibitem{OM08} 
	\textnormal{R. E. O'Malley}, 
	\textit{Singularly perturbed linear two-point boundary value problems.} 
	SIAM Rev. {\bf 50} (2008), no. 3, pp 459-482.

\bibitem{RB19}
Nasim Rahaman, Aristide Baratin, Devansh Arpit, Felix Draxler, Min Lin, Fred Hamprecht, Yoshua Bengio, and Aaron Courville.
On the Spectral Bias of Neural Networks.  
{\em Proceedings of the 36th International Conference on Machine Learning}, 97 (2019). 

\bibitem{RK18}
M. Raissi and G. E. Karniadakis.  
\newblock {Hidden physics models: machine learning of nonlinear partial didderential equations}.
{\em J. Comput. Phys.}, 
357 (2018), pp. 125--141.

\bibitem{RPK19}
M. Raissi, P. Perdikaris, and G. E. Karniadakis.   
\newblock {Physics-informed neural networks: a deep learning framework for solving forward and inverse problems involving nonlinear partial differential equations}.  
{\em J. Comput. Phys.}, 378 (2019),
pp. 686--707.


\bibitem{RBPK17}
S. H. Rudy, S. L. Brunton, J. L. Proctor, J. N. Kutz.  
Data-driven discovery of partial differential equations.  
{\em Science advances}, 3 (4) (2017).


\bibitem{SDBSK20}
K. Shukla, P. C. Di Leoni, J. Blackshire, D. Sparkman, G. E. Karniadakis. 
Physics-informed neural network for ultrasound nondestructive quantification of surface breaking cracks.  
{\em Journal of Nondestructive Evaluation}, 39 (3) (2020) 1--20.

	
\bibitem{SK87} 
	\textnormal{S. Shih and R. B. Kellogg}, 
	\textit{Asymptotic analysis of a singular perturbation problem.} 
	SIAM J. Math. Anal. {\bf 18} (1987), pp . 1467-1511.



\bibitem{SS18}
J. Sirignano, K. Spiliopoulos.  
Dgm: A deep learning algorithm for solving partial differential equations.  
{\em Journal of computational physics}, 375 (2018) 1339--1364.


\bibitem{FSF22}
Mario De Florio, Enrico Schiassi, and Roberto Furfaro.  
\newblock {Physics-informed neural networks and functional interpolation for stiff chemical kinetics}. 
{\em Chaos} 32, 063107 (2022).


\bibitem{QCJX21}
T. Qin, Z. Chen, J. D. Jakeman, D. Xiu.  
Data-driven learning of nonautonomous systems.  
{\em SIAM Journal on Scientific Computing}, 43 (3) (2021).


\bibitem{WSP22}
S. Wang, S. Sankaran, P. Perdikaris.  
Respecting causality is all you need for training physics-informed neural networks.  
{\em arXiv preprint}, arXiv:2203.07404 (2022).


\bibitem{XZRW21}
R. Xu, D. Zhang, M. Rong, N. Wang.  
Weak form theory-guided neural network (tgnn-wf) for deep learning of subsurface single-and two-phase flow.  
{\em Journal of Computational Physics}, 436 (2021).



\bibitem{YMK21}
L. Yang, X. Meng, G. E. Karniadakis.  
B-pinns: Bayesian physics-informed neural networks for forward and inverse pde problems with noisy data.  
{\em Journal of Computational Physics}, 425 (2021).


\bibitem{Yu18}
B. Yu, et al.  
The deep Ritz method: a deep learning-based numerical algorithm for solving variational problems.  
{\em Communications in Mathematics and Statistics}, 6 (1) (2018).










\end{thebibliography}
\end{document}